%% file: voltpWH.tex
\documentclass[11pt,a4wide]{article}

\usepackage{geometry}
\usepackage{bm} % pdreesen \bm math boldface
\usepackage{graphicx}

\usepackage{amssymb} % pdreesen amssymb for mathbb
\usepackage{amsmath}

\usepackage{tikz}
\usetikzlibrary{calc,arrows}

\renewcommand{\vec}[1]{\mathbf{#1}}

\newcommand{\mat}[1]{\mathbf{#1}}
\DeclareMathAlphabet\mathbfcal{OMS}{cmsy}{b}{n} % boldface caligraphic math
\newcommand{\ten}[1]{\mathbfcal{#1}}
\usepackage{relsize} % pdreesen for \mathsmaller
\renewcommand{\top}{\mathsmaller{T}}

\newcommand{\llbrack}{\left[\!\left[}
\newcommand{\rrbrack}{\right]\!\right]}
\newcommand{\llbrackbig}{\Big[\!\!\Big[}
\newcommand{\rrbrackbig}{\Big]\!\!\Big]}

\newcommand{\citep}[1]{\cite{#1}}%{\color{red}(#1)}}
\newcommand{\citet}[1]{\cite{#1}}%{\color{blue}(#1)}}

\everymath{\displaystyle}

\begin{document}

\input{pretikz.tex}

\title{Modeling Parallel Wiener-Hammerstein Systems\\Using Tensor Decomposition of Volterra Kernels}

%\titlerunning{Modeling Parallel Wiener-Hammerstein Systems Using Tensors}

%\toctitle{Modeling Parallel Wiener-Hammerstein Systems Using Tensor Decomposition of Volterra  Kernels}

\author{Philippe~Dreesen$^1$, David~Westwick$^2$, Johan~Schoukens$^1$, Mariya~Ishteva$^1$\\ \\ $^1$ Vrije Universiteit Brussel, Dept~VUB-ELEC, Brussels, Belgium\\ $^2$ University of Calgary, Dept~Electrical and Computer Engineering, Calgary, Canada\\ {\tt{philippe.dreesen@gmail.com}}}
%\authorrunning{Philippe Dreesen et al.}
%\tocauthor{Philippe Dreesen, David Westwick, Johan Schoukens and Mariya Ishteva}

\maketitle

\begin{abstract}
Providing flexibility and user-interpretability in nonlinear system identification can be achieved by means of block-oriented methods. 
One of such block-oriented system structures is the parallel Wiener-Hammerstein system, which is a sum of Wiener-Hammerstein branches, consisting of static nonlinearities sandwiched between linear dynamical blocks. 
Parallel Wiener-Hammerstein models have more descriptive power than their single-branch counterparts, but their identification is a non-trivial task that requires tailored system identification methods. 
In this work, we will tackle the identification problem by performing a tensor decomposition of the Volterra kernels obtained from the nonlinear system. 
We illustrate how the parallel Wiener-Hammerstein block-structure gives rise to a joint tensor decomposition of the Volterra kernels with block-circulant structured factors. 
The combination of Volterra kernels and tensor methods is a fruitful way to tackle the parallel Wiener-Hammerstein system identification task. 
In simulation experiments, we were able to reconstruct very accurately the underlying blocks under noisy conditions. 

\end{abstract}

\section{Introduction}\label{sec:intro}
System identification is the art of building dynamical models from noisy measurements of input and output data.
Linear system identification is a well-established discipline \cite{ljung1999sitftu,pintelon2012siafda,vanoverschee1996siflstia} and has yielded successful applications in a wide variety of fields. 
In the last decades, the use of nonlinear models has become more important in order to capture the nonlinear effects of the real world. 
Many different nonlinear identification methods have been proposed~\cite{giannakis2001abonsi,sjoberg1995nbmisiauo}, but very often these solutions are either tailored to a specific application, or are too complex to understand or study. 

In this paper, we will draw ideas from two nonlinear system identification approaches and try to combine the benefits of both. 
The first approach tackles the disadvantage of increased complexity of nonlinear models by considering block-oriented models~\cite{giri2010bnsi}, which combine flexibility with user-interpretation by interconnecting linear dynamical blocks and static nonlinear functions. 
Unfortunately, even simple block-oriented models, such as Wiener (cascade of linear-nonlinear), Hammerstein (cascade nonlinear-linear) or Wiener-Hammerstein (cascade linear-nonlinear-linear) require an iterative optimization on a non-convex objective function, and identification procedures that are tailored towards a specific block-structure.
The second approach that we will use is the Volterra model \cite{boyd1985volterra,schetzen1980}, which is an extension of the well-known impulse response model for linear dynamical systems. 
Volterra models take into account higher-order polynomial nonlinearities and can thus be seen as a generalization of the Taylor series expansion for nonlinear dynamical systems.
The advantages of Volterra models are that any fading-memory system can be approximated to an arbitrary degree of accuracy \cite{boyd1985volterra,schoukens2005} and the parameter estimation task is a linear problem. 
The disadvantages are that the resulting models contain a very high number of parameters and thus cannot be given physical interpretation.

We generalize the earlier results of \cite{degoulart2016,favier2006pt1} on Wiener-Hammerstein system identification using tensor decompositions in two ways:
First, we show that the case of parallel branches for a fixed degree $d$ gives rise to a canonical polyadic decomposition with block-structured factors. 
The study of parallel Wiener-Hammerstein systems is useful, as they are universal approximators, where as single-branch Wiener-Hammerstein systems are not \cite{palm1979,schetzen1980}. 
Second, we jointly consider Volterra kernels of several degrees by means of the structured data fusion framework of \cite{sorber2015sdf} and solve the problem as a joint structured tensor decomposition. 
By simultaneously decomposing several Volterra kernels, the available information is used maximally. 
The presented method is implemented by means of structured data fusion \cite{sorber2015sdf} using Tensorlab 3.0 \cite{tensorlab3}, and is validated on simulation experiments.

The paper is organized as follows. 
In Section~\ref{sec:volterra} we illustrate the link between the Volterra kernels and tensors and introduce the canonical polyadic decomposition for tensors. 
Section~\ref{sec:pWHstruct} illustrates how the Volterra kernels of a parallel Wiener-Hammerstein system have a natural connection to the canonical polyadic decomposition with block-circulant factors.
This ultimately leads to a joint structured canonical polyadic decomposition of the Volterra kernels of several orders that solve the parallel Wiener-Hammerstein identification task. 
We validate the method on numerical simulation examples in Section~\ref{sec:exp}. 
In Section~\ref{sec:conc} we draw the conclusions.

\section{Volterra kernels, tensors and tensor decomposition}\label{sec:volterra}
\subsection{The Volterra model for nonlinear systems}
We consider single-input-single-output systems that map an input $u$ at time instance $k$ onto an output $y$ at time instance $k$. 
The Volterra series expansion generalizes the well-known convolution operator for linear dynamical systems to the nonlinear case.
Essentially the Volterra model of a nonlinear system expresses the output $y(k)$ as a polynomial function of time-shifted input variables $u(k), u(k-1), \ldots, u(k-m)$, with $m$ denoting the memory length of the model.
Formally we can write the Volterra model as
\begin{equation}
y(k) = \sum_{d=1}^D \left( \sum_{s_1,\ldots,s_d=0}^m H_d(s_1, \ldots, s_d) u(k-s_1) \cdots u(k- s_d)  \right),
\label{eq:ykVolterra}
\end{equation}
where $H_d(\cdot,\ldots,\cdot)$ denotes Volterra kernel of degree $d$. 
The Volterra series expansion allows for representing a large class of nonlinear systems up to an arbitrary degree of accuracy \cite{boyd1985volterra}.

\subsection{From polynomials to tensors}
It is well-known that multivariate homogeneous polynomials can be identified to symmetric higher-order tensors \cite{comon2008symtensors}.
For instance, we may represent the quadratic polynomial $p(x_1,x_2) = 5 x_1^2 -8 x_1 x_2 + x_2^2$ as a matrix multiplied from both sides by a vector containing $x_1$ and $x_2$ as 
$$ \begin{array}{rcl} p(x_1,x_2) &=&  5 x_1^2 -8 x_1 x_2 + x_2^2 \\ \\
    &=& \left[ \begin{array}{cc} x_1 & x_2 \end{array} \right]
        \left[ \begin{array}{rrr} 5 && -4\\-4 && 1 \end{array} \right]
\left[ \begin{array}{c} x_1 \\ x_2 \end{array} \right]. \end{array} 
$$
In general, we may thus write a (nonhomogeneous) polynomial as 
\begin{equation} 
    p(x_1, \ldots, x_n) = p_0 + \vec{x}^\top \vec{p}_1 + \vec{x}^\top \mat{P}_2 \vec{x} + \ten{P}_3 \times_1 \vec{x}^\top \times_2 \vec{x}^\top \times_3 \vec{x}^\top + % \ten{H}_4 \times_1 \vec{u}^\top \times_2 \vec{u}^\top \times_3 \vec{u}^\top \times_4 \vec{u}^\top + 
%    y(k) = \vec{u}^\top \vec{h}_1 + \vec{u}^\top \mat{H}_2 \vec{u} + \ten{H}_3 \times_1 \vec{u}^\top \times_2 \vec{u}^\top \times_3 \vec{u}^\top + % \ten{H}_4 \times_1 \vec{u}^\top \times_2 \vec{u}^\top \times_3 \vec{u}^\top \times_4 \vec{u}^\top + 
    \ldots,
    \label{eq:polyistensor}
\end{equation}
    where $\vec{x} = \left[ \begin{array}{ccc} x_1 & \ldots & x_n \end{array} \right]^\top$ and $\times_n$ is the $n$-mode product defined as follows. 
Let $\ten{X}$ be a $I_1 \times I_2 \times \cdots \times I_N$ tensor, 
and let $\vec{u}^\top$ be an $1 \times I_n$ row vector, 
then we have $$ \left( \ten{X} \times_n \vec{u}^\top \right)_{i_1 \cdots i_{n-1} i_{n+1} \cdots i_N } = \sum_{i_n=1}^{I_n} x_{i_1 i_2 \cdots i_N} u_{i_n}.$$
Notice that the result is a tensor of order $N-1$, as mode $n$ is summed out. 

\subsection{Canonical Polyadic Decomposition}
It is often useful to decompose a tensor into simpler components, and for the proposed method we will use the canonical polyadic decomposition. 
The canonical polyadic decomposition \cite{carroll1970,harshman1970,kolda2009tdaa} (also called CanDecomp or PARAFAC) expresses the tensor $\ten{T}$ as a sum of rank-one terms as
$$
    \ten{T} = \sum_{i=1}^R \vec{a}_r \circ \vec{b}_r \circ \vec{c}_r,
$$
with $\circ$ denoting the outer product and the number of components $R$ denoting the CP rank of the tensor $\ten{T}$.  
Often a short-hand notation is used as $$\ten{T} = \llbrack \mat{A}, \mat{B}, \mat{C} \rrbrack,$$
where $\mat{A} = \left[ \begin{array}{ccc} \vec{a}_1 & \cdots & \vec{a}_R \end{array}\right]$ and $\mat{B}$ and $\mat{C}$ are defined likewise.

%\newpage
\section{Parallel Wiener-Hammerstein as tensor decomposition}\label{sec:pWHstruct}
For self-containment we will first review and rephrase a result of \cite{favier2006pt1} that connects the Volterra kernel of a Wiener-Hammerstein system to a canonical polyadic decomposition with circulant-structured factor matrices. 
Afterwards, we will generalize this to the parallel case and then we show that the entire problem leads to a joint and structured canonical polyadic decomposition.

\subsection{Wiener-Hammerstein as structured tensor decomposition} 
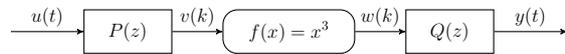
\begin{figure}[h!]
\begin{center}
\tikz \node[scale=0.5]{ \input{fig-WH.tex} };
 \end{center}
	\vspace*{-0.2cm}
 \caption{A Wiener-Hammerstein system with input signal $u(k)$ and output signal $y(k)$ contains  a static nonlinear function $f(\cdot)$ that is sandwiched between the FIR filters $P(z)$ and $Q(z)$.}
 \label{fig:WH}
\end{figure}

Let us try to understand how the canonical polyadic decomposition shows up in modeling a Wiener-Hammerstein system.
Consider a (single-branch) Wiener-Hammerstein system as in Figure~\ref{fig:WH} with FIR filters $P(z)$ and $Q(z)$ with memory lengths $m_P$ and $m_Q$, respectively, and a static nonlinearity $f (x) = x^3$. 
The output $y(k)$ of the Wiener-Hammerstein model is obtained by passing the signal $w(k)$ through the filter $Q(z)$. 
We can write this as 
\begin{equation}
\begin{array}{rcl} 
y (k) 
      &=& \left[ \begin{array}{ccc} w(k) & \cdots & w(k-m_Q) \end{array} \right] \left[ \begin{array}{cccc} 1 & q_1 & \cdots & q_{m_Q} \end{array} \right]^\top \\ \\ 
&=& \vec{w}^\top \vec{q}, \\ \\  
\end{array}
\label{eq:ywTq}
\end{equation}
where we fixed the first filter coefficient $q_0 = 1$ in order to ensure uniqueness of the identified model.
The signal $w(k)$ is given by the expression $w(k) = f(v(k))$, or in this case $w(k) = v^3(k)$.
To obtain $v(k), \ldots, v(k-m_Q)$ from $u(k)$, we require of $u(k)$ the samples $k$ down to $k-m_Q-m_P$. 
This convolution operation can be expressed as a matrix equation as 
$$
\begin{array}{rcl}
    %\vec{v}^\top &=& \vec{u}^\top \mat{P} \\ \\
\left[ \begin{array}{ccc} v(k) & \cdots & v(k-m_Q) \end{array} \right] 
    &=&
\left[ \begin{array}{ccc} u(k) & \cdots & u(k-m_Q-m_P) \end{array} \right] 
    \left[ \begin{array}{ccc} 1 & & \\ p_1 & \ddots & \\ \vdots & \ddots & 1 \\ p_{m_P} & & p_1 \\ & \ddots & \vdots \\ &  & p_{m_P} \end{array} \right]\\ \\
\vec{v}^\top&=& \vec{u}^\top\mat{P},\\
\end{array}
$$ 
with the circulant matrix $\mat{P}$ of size $m_P + m_Q + 1 \times m_Q + 1$. 
Notice that we fixed the first coefficient $p_0= 1$ for uniqueness purposes.
The matrix $\mat{P}$ will turn out to play a central role in the canonical polyadic decomposition of the Volterra kernels of a Wiener-Hammerstein system.

By fixing both $q_0=1$ and $p_0=1$, we are excluding the possiblity that there is a pure delay present in the system. 
The presence of a delay in the system can be accounted for by setting $p_0=1$ and then performing a scaling on the nonlinearity, rather than on $q_0$.
In case the system has a delay this will lead to an estimated $q_0 \approx 0$.  
However, for notational convenience we have chosen $p_0=q_0=1$ in the remainder of this paper, but a more general scaling strategy to ensure uniqueness is possible and compatible with the presented method.

For the current Wiener-Hammerstein system we have $f(x)=x^3$, and hence $y (k) = \ten{H}_3 \times_1 \vec{u}^\top\times_2 \vec{u}^\top\times_3 \vec{u}^\top$. 
In \cite{favier2006pt1} it is shown that the Volterra kernel can be written as the canonical polyadic decomposition $ \ten{H} = \llbrack \mat{P}, \mat{P}, \mat{P} \operatorname{diag} (\vec{q}) \rrbrack $, which we can also write in a more symmetrical expression by extracting $\vec{q}^\top$ into an extra mode as 
\begin{equation}\ten{H} = \llbrack \mat{P}, \mat{P}, \mat{P}, \vec{q}^\top \rrbrack.\label{eq:cpdWH}\end{equation}
%of which the dimension is given by $(m_P + m_Q + 1) \times (m_P + m_Q + 1) \times (m_P + m_Q + 1)  \times 1 = (m_P + m_Q + 1) \times (m_P + m_Q + 1) \times (m_P + m_Q + 1)$.   
This fact can be appreciated by considering output $y(k)$ as  
$$\begin{array}{rcl} y(k) &=& \ten{H} \times_1 \vec{u}^\top \times_2 \vec{u}^\top \times_3 \vec{u}^\top \\ \\ &=& \llbrack \vec{u}^\top \mat{P}, \vec{u}^\top \mat{P}, \vec{u}^\top \mat{P}, \vec{q}^\top \rrbrack \\ \\ &=& \llbrack \vec{v}^\top, \vec{v}^\top, \vec{v}^\top, \vec{q}^\top \rrbrack \\ \\ &=& \sum_{i=0}^{m_Q} q(i) v^3(k-i), \end{array}$$ 
in which we recognize the convolution of the impulse response of $Q(z)$ with the time-shifted samples $v^3 (k)$ as in (\ref{eq:ywTq}).  

In case of a general polynomial function $f(x)$, the same reasoning can be developed for each degree $d$, which will lead to a structured canonical polyadic decomposition of the degree-$d$ Volterra kernel as in (\ref{eq:cpdWH}). 
For instance, if $f(x)=a x^2 + b x^3$, we find the following expressions
$$ \begin{array}{rcl} \mat{H}_2 &=& a \llbrack \mat{P}, \mat{P}, \vec{q}^\top \rrbrack, \\ \\ 
\ten{H}_3 &=& b \llbrack \mat{P}, \mat{P}, \mat{P}, \vec{q}^\top \rrbrack. \end{array} $$
In Section~\ref{sec:jointdecomp} we will discuss how this leads to a joint tensor decomposition.  

\subsection{Parallel Wiener-Hammerstein structure}
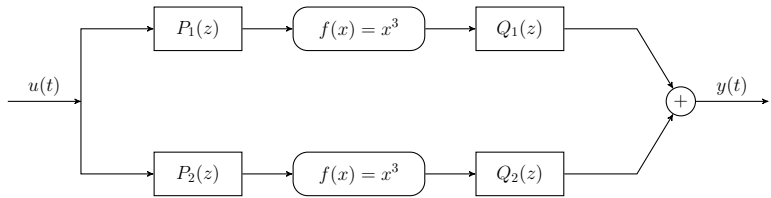
\begin{figure}[b]
    \begin{center}
        \tikz \node [scale=0.5]{\input{fig-pWH.tex}};
    \end{center}
	\vspace*{-0.2cm}
    \caption{An example of a two-branch parallel Wiener-Hammerstein system having an identical nonlinear function $f_1(x)=f_2(x)=x^3$.}
    \label{fig:pWH}
\end{figure}

To understand how we can extend these results to the parallel case, let us consider a two-branch parallel Wiener-Hammerstein system where both branches have an identical nonlinearity $f_1(x)=f_2(x)=x^3$, as in Figure~\ref{fig:pWH}.
To avoid a notational overload, we will assume for the remainder of this paper that the memory lengths of all filters $P_i(z)$ are $m_P$, and likewise for the filters $Q_i(z)$ the lenghts are $m_Q$. 
The summation of the two branches leads to 
$$
\begin{array}{rcl}
    \ten{H}_3 &=& \llbrack \mat{P}_1, \mat{P}_1, \mat{P}_1, \vec{q}_1^\top \rrbrack + \llbrack \mat{P}_2, \mat{P}_2, \mat{P}_2, \vec{q}_2^\top \rrbrack \\ \\
         &=& \llbrackbig \left[ \begin{array}{cc} \mat{P}_1 & \mat{P}_2 \end{array} \right], \left[ \begin{array}{cc} \mat{P}_1 & \mat{P}_2 \end{array} \right],  \left[ \begin{array}{cc} \mat{P}_1 & \mat{P}_2 \end{array} \right], \left[ \begin{array}{cc} \vec{q}_1^\top & \vec{q}_2^\top\end{array} \right] \rrbrackbig,\end{array}
$$ 
with $\mat{P}_i$ and $\vec{q}_i$ defined similar as in the single-branch case. 
We may include a scaling of branch one by a scalar $c_1$ (i.e., $f_1(x) = c_1 x^3$ and branch two by a scalar $c_2$ (i.e., $f_2(x) = c_2 x^3$ by introducing an additional mode as
\begin{equation}
    \ten{H}_3 = \llbrackbig \left[ \begin{array}{cc} \mat{P}_1 & \mat{P}_2 \end{array} \right], \left[ \begin{array}{cc} \mat{P}_1 & \mat{P}_2 \end{array} \right],  \left[ \begin{array}{cc} \mat{P}_1 & \mat{P}_2 \end{array} \right], \left[ \begin{array}{cc} \vec{q}_1^\top & \vec{q}_2^\top\end{array} \right], \left[ \begin{array}{cc} c_1 \vec{1}_{m_Q+1}^\top & c_2 \vec{1}_{m_Q+1}^\top \end{array}\right] \rrbrackbig.
\label{eq:pWHCPDd3}
\end{equation}
%where `factor' $\left[ \begin{array}{cc} c_1 \vec{1}_{m_Q+1}^\top & c_2 \vec{1}_{m_Q+1}^\top \end{array} \right]$ corresponds to the multiplication of all $\mat{P}_1$ and $\vec{q}_1$ with $c_1$, and all $\mat{P}_2$ and $\vec{q}_2$ with $c_2$, respectively. 
Introducing the extra factor $\left[ \begin{array}{cc} \vec{1}_{m_Q+1}^\top c_1 & \vec{1}_{m_Q+1}^\top c_2 \end{array}\right]$ does not change the size of the tensor, since it introduces a mode with dimension one. 
Formally, if we let $m=m_P + m_Q + 1$ denote the memory length of the Volterra model, the tensor in (\ref{eq:pWHCPDd3}) has size $m \times m \times m \times 1 \times 1$ which is equivalent to $m \times m \times m$. 
        
\subsection{Coupled tensor and matrix decompositions}\label{sec:jointdecomp}
%The general expression for the output $y (k)$ can be written as 
%$$ y(k) = \vec{u}^\top \vec{h}_1 + \vec{u}^\top \mat{H}_2 \vec{u} + \ten{H}_3 \times_1 \vec{u}^\top \times_2 \vec{u}^\top \times_3 \vec{u}^\top + \ldots.$$
The Volterra kernels of the parallel Wiener-Hammerstein model for a particular order $d$ can be decomposed as a structured canonical polyadic decomposition.
Hence, if the Volterra kernels of multiple orders are available, a joint decomposition of multiple Volterra kernels should be performed. 

Ultimately, we find that the $R$-branch parallel Wiener-Hammerstein identification task is solved by minimizing the cost criterion
\begin{equation}
    \begin{array}{l}
\operatorname*{minimize}_{\mat{P}, \vec{q}, \vec{c}} 
\left\| \vec{h}_1 - \llbrack \mat{P}, \vec{q}^\top, \vec{c}_1^\top \rrbrack \right\|_2^2 
+
\left\| \mat{H}_2 - \llbrack \mat{P}, \mat{P}, \vec{q}^\top, \vec{c}_2^\top \rrbrack \right\|_F^2 \\ \\
\qquad \qquad  \qquad  \qquad  \qquad  \qquad \qquad + \left\| \ten{H}_3 - \llbrack \mat{P}, \mat{P}, \mat{P}, \vec{q}^\top, \vec{c}_3^\top \rrbrack \right\|_F^2 + \ldots, \end{array}
\label{eq:tenmat}
\end{equation}
where 
$$
\begin{array}{rcl}
    \mat{P} &=& \left[ \begin{array}{ccc} \mat{P}_1 & \cdots & \mat{P}_R \end{array} \right],\\ \\
        \vec{q}^\top &=& \left[ \begin{array}{ccc} \vec{q}_1^\top & \cdots & \vec{q}_R^\top \end{array} \right], \\ \\
            \vec{c}_d^\top &=& \left[ \begin{array}{ccc} c_{1d} \vec{1}_{m_Q+1}^\top & \cdots & c_{Rd} \vec{1}_{m_Q+1}^\top \end{array} \right]. \\ \\
            %\vec{c}_3^\top &=& \left[ \begin{array}{cc} c_{13} \vec{1}_{m_Q+1}^\top & c_{23} \vec{1}_{m_Q+1} \end{array} \right]. 
\end{array}
$$
The factor matrices $\mat{P}$ and $\vec{q}^\top$ are shared among all joint decompositions while the constants $\vec{c}_d^\top$ depend on the order $d$ of the considered Volterra kernel. 
Joint and structured factorizations like (\ref{eq:tenmat}) can be solved in the framework of structured data fusion \cite{sorber2015sdf}. 

\section{Numerical results}\label{sec:exp}
In this section we validate the proposed identification method on a simulation example. 
Numerical experiments were performed using MATLAB and structured data fusion \cite{sorber2015sdf} in Tensorlab 3.0 \cite{tensorlab3} (MATLAB code is available on request). 

\begin{figure}
    \begin{center}
        \tikz \node[scale=0.5]{\input{fig-pWHnoise.tex}};
    \end{center}
	\vspace*{-0.2cm}
    \caption{A two-branch parallel Wiener-Hammerstein system with output noise.}
    \label{fig:pWHnoise}
\end{figure}
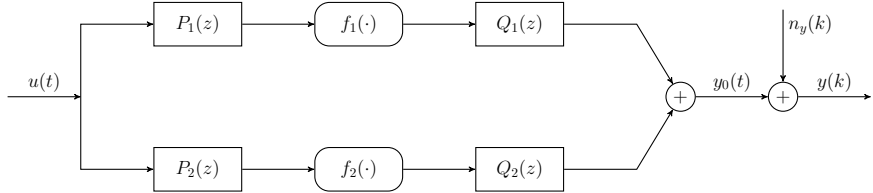

We consider a parallel Wiener-Hammerstein system having two branches with second and third-degree polynomial nonlinearities (Figure~\ref{fig:pWHnoise}).
The finite impulse response coefficients of $P_i(z)$ and $Q_i(z)$ are chosen as sums of decreasing exponentials with lengths $m_P = m_Q = 10$.
The coefficients $c_{i2}$ and $c_{i3}$ are drawn from a normal distribution $N(0,0.1^2)$. 
The input signal is a Gaussian white noise sequence ($u(k) \sim N(0,0.7^2)$) and applied without noise to the system, for $k=1,\ldots,10,000$. 
The outputs $y(k) = y_0(k) + n_y(k)$ are contaminated with additive Gaussian noise $n_y$ with a signal-to-noise ratio of 10 dB.

The Volterra kernels of degrees two and three are estimated using a standard least-squares method on the basis of the inputs $u(k)$ and the noisy outputs $y(k)$. 
The memory lengths $m_P=m_Q=10$ give rise to $21 \times 21$ second-order Volterra kernel and a $21 \times 21 \times 21$ third-order Volterra kernel, having in total $231 + 1771 = 2002$ kernel elements.

\begin{figure}[!htb]
\begin{minipage}[c]{0.5\textwidth}
\begin{center}
\includegraphics[width=1\textwidth]{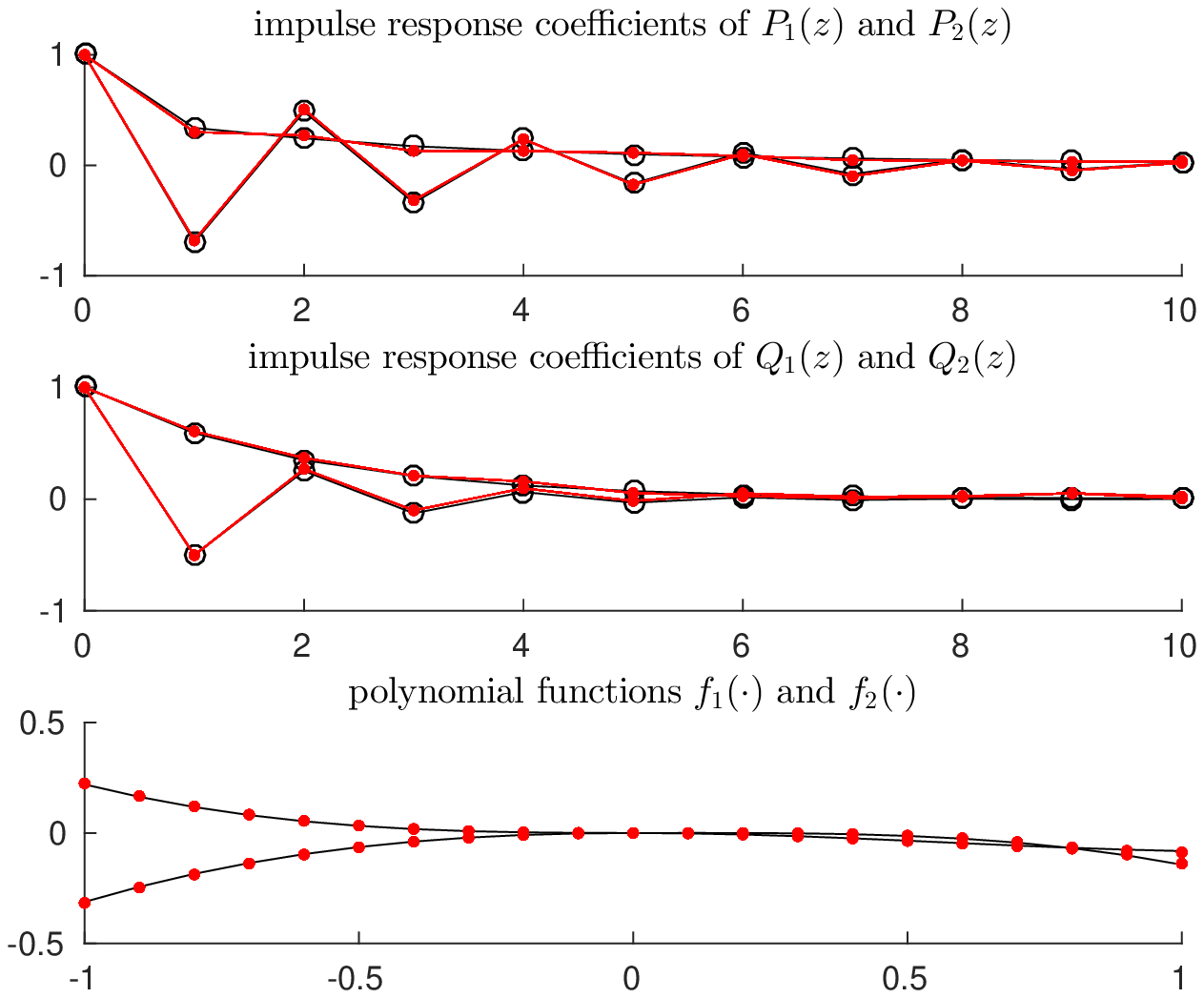}
(a)
\end{center}
\end{minipage}
\begin{minipage}[c]{0.5\textwidth}
\begin{center}
\includegraphics[width=1\textwidth]{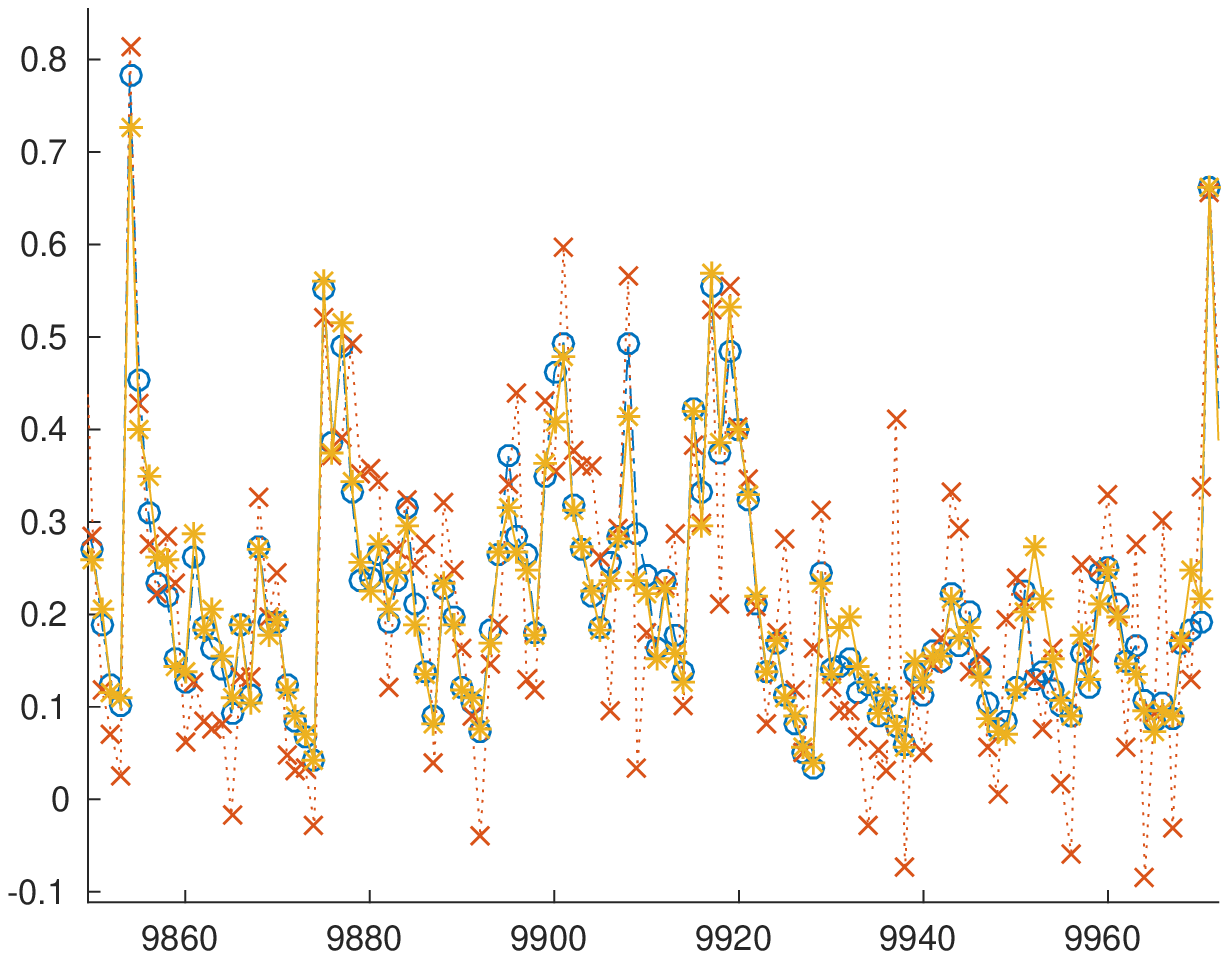}
(b)
\end{center}
\end{minipage}
	\vspace*{-0.2cm}
\caption{The method succeeds in retrieving the parallel Wiener-Hammerstein system parameters about 10\% of the cases. Typical results of successful completion are shown in (a). The output of the identified parallel Wiener-Hammerstein system ({\color{orange}$*$}) reconstructs very accurately the noiseless output signal $y_0$ ({\color{blue}$\circ$}), while the Volterra kernels were computed from the noisy data with 10 dB SNR ({\color{red}$\times$}).}
\label{fig:ex}
\end{figure}

The joint matrix and tensor decomposition with structured factors is then performed, returning the parameters of the parallel Wiener-Hammerstein system.
We have performed a Monte Carlo experiment with 100 re-initializations of the optimization routine and were able to retrieve the true underlying system parameters in about 10~\% out of the time. 
In Figure~\ref{fig:ex}~(a) we show a typical result of a successful identification that was able to retrieve the underlying system parameters accurately.  
The reconstructed outputs of the identified model together with the true outputs and the noisy measurements from which the Volterra kernels were estimated are shown in Figure~\ref{fig:ex}~(b).  

In other experiments we have observed that the success rate of the algorithm increases for lower noise levels and shorter filter lengths. 
If the method failed consistently on a system, this could almost always be understood from problem-specific system properties, such as similar impulse responses in the filters, dominance of one of the branches (i.e., unbalanced values of the coefficients of $f_i$), etc. 
Nevertheless, for extending the method to larger-scale problems, the issue of having good initial estimates becomes more relevant. 
Possibly the work of \cite{frederik} which focuses on block-circulant structured canonical polyadic decomposition may provide good starting points for initializing our method.

Finally we would like to remark that in simulations, one can easily compare the retrieved system parameters (or simulated output data) with the underlying true system parameters or signals, which is obviously not possible in real experiments. 
However, it is worth mentioning that the tensor approximation error~(\ref{eq:tenmat}) was strongly correlated with the error between the simulated output and noiseless output data, which provides a proxy to select the best model among a set of candidates. 

\section{Conclusions}\label{sec:conc}
The joint decomposition of Volterra kernels with structured factors is able to retrieve up to a high degree of accuracy the underlying system parameters. 
The success rate of the method decreases as the noise level and number of parameters in the system grows, but even up to moderately long impulse responses (memory length of ten samples), the method was successful in about 10 \% of the re-initializations of the optimization routine. 
An important challenge for future work is to obtain good initializations, as this becomes an important issue when considering filters with longer memories and/or higher noise levels.

\section*{Acknowledgments}\label{sec:ack}
This work was supported in part by the Fund for Scientific Research (FWO-Vlaanderen), by the Flemish Government (Methusalem), the Belgian Government through the Inter-university Poles of Attraction (IAP VII) Program, by the ERC Advanced Grant SNLSID under contract 320378 and by FWO project G.0280.15N.
The authors want to thank Otto Debals and Nico Vervliet for technical assistance for the use of Tensorlab/SDF and the suggestion to extract the vector $\vec{q}$ into an additional tensor mode.

\bibliographystyle{plain}

\bibliography{references}

\end{document}

%% file: pretikz.tex
% this goes AFTER \begin{document}

%\begin{document}

% Define lengths
\newlength{\wLTI}
\newlength{\wSNL}
\newlength{\wBigSNL}
\newlength{\wModel}
\newlength{\wShortNamedSignal}
\newlength{\wUnnamedSignal}
\newlength{\hLTI}
\newlength{\dyLTI}
\newlength{\sSum}

% Set lengths
\setlength{\wLTI}{6em}
\setlength{\wSNL}{\wLTI}
\setlength{\wBigSNL}{9em}
\setlength{\wModel}{\wLTI}
\setlength{\wShortNamedSignal}{5em}
\setlength{\wUnnamedSignal}{2em}
\setlength{\hLTI}{3em}
\setlength{\dyLTI}{2em}
\setlength{\sSum}{2em}

% Define colors
\colorlet{LTIcolor}{black}
\colorlet{SNLcolor}{black}
\colorlet{MODELcolor}{black}

% Define block styles
\tikzstyle{SISO_LTI} 	= [draw, fill=none, text centered, minimum height=\hLTI, minimum width=\wLTI, LTIcolor]
\tikzstyle{SISO_SNL}	= [draw, fill=none, text centered, minimum height=\hLTI, minimum width=\wSNL, rounded corners=1em, SNLcolor]
\tikzstyle{SISO_BigSNL}	= [draw, fill=none, text centered, minimum height=\hLTI, minimum width=\wBigSNL, rounded corners=1em, SNLcolor]
\tikzstyle{sum} 	= [draw, fill=none, text centered, minimum size=\sSum, circle]

%% file: fig-WH.tex
% wiener-hammerstein figure with intermediate signals v and w 
\begin{tikzpicture}[font=\Large, >=stealth', thick, auto]
	% general coordinates
	\node (u_start)	[coordinate]{};
	\node (u_end)		[coordinate, right of=u_start, node distance=\wShortNamedSignal]{};
	\node (u1_start)	[coordinate, above of=u_end, node distance=(\hLTI+\dyLTI)]{};
	\node (u3_start)	[coordinate, below of=u_end, node distance=(\hLTI+\dyLTI)]{};
	\node (front_LTI1_c)	[coordinate, right of=u1_start, node distance=(\wShortNamedSignal+\wLTI/2)]{};
	\node (front_LTI2_c)	[coordinate, right of=u_end, node distance=(\wShortNamedSignal+\wLTI/2)]{};
	\node (front_LTI3_c)	[coordinate, right of=u3_start, node distance=(\wShortNamedSignal+\wLTI/2)]{};
	\node (SNL1_c)	[coordinate, right of=front_LTI1_c, node distance=(\wLTI/2+\wShortNamedSignal+\wSNL/2)]{};
	\node (SNL2_c)	[coordinate, right of=front_LTI2_c, node distance=(\wLTI/2+\wShortNamedSignal+\wSNL/2)]{};
	\node (SNL3_c)	[coordinate, right of=front_LTI3_c, node distance=(\wLTI/2+\wShortNamedSignal+\wSNL/2)]{};
	\node (back_LTI1_c)	[coordinate, right of=SNL1_c, node distance=(\wSNL/2+\wShortNamedSignal+\wLTI/2)]{};
	\node (back_LTI2_c)	[coordinate, right of=SNL2_c, node distance=(\wSNL/2+\wShortNamedSignal+\wLTI/2)]{};
	\node (back_LTI3_c)	[coordinate, right of=SNL3_c, node distance=(\wSNL/2+\wShortNamedSignal+\wLTI/2)]{};
	\node (y1_end)	[coordinate, right of=back_LTI1_c, node distance=(\wLTI/2+\wShortNamedSignal)]{};
	\node (y2_end)	[coordinate, right of=back_LTI2_c, node distance=(\wLTI/2+\wShortNamedSignal)]{};
	\node (y3_end)	[coordinate, right of=back_LTI3_c, node distance=(\wLTI/2+\wShortNamedSignal)]{};
	\node (sum_c)		[coordinate, right of=y2_end, node distance=(\wUnnamedSignal+\sSum/2)]{};
	\node (y_end)		[coordinate, right of=sum_c, node distance=(\sSum/2+\wShortNamedSignal)]{};

	% specific nodes
        \node (LTI1)		at (front_LTI1_c) 	[SISO_LTI]{$P(z)$};
        \node (SNL)		at (SNL1_c)		[SISO_BigSNL]{$f(x)=x^3$};
        \node (LTI2)		at (back_LTI1_c)	[SISO_LTI]{$Q(z)$};
		
	% signals
	\path[-]	(u1_start)	-- node {$u(t)$}	(LTI1);
	\path[-]	(LTI2)		-- node {$y(t)$}	(y1_end);

	% arrows
        \draw[->] 	(u1_start)	-- (LTI1);
        \draw[->] 	(LTI1)		-- node {$v(k)$} 	(SNL);
        \draw[->]	(SNL)		-- node {$w(k)$} 	(LTI2);
	\draw[->]	(LTI2)		-- 	(y1_end);
\end{tikzpicture}

%% file: fig-pWH.tex
% parallel wiener-hammerstein figure 
\begin{tikzpicture}[font=\Large, >=stealth', thick, auto]
	% general coordinates
	\node (u_start)	[coordinate]{};
	\node (u_end)		[coordinate, right of=u_start, node distance=\wShortNamedSignal]{};
	\node (u1_start)	[coordinate, above of=u_end, node distance=(\hLTI+\dyLTI)]{};
	\node (u3_start)	[coordinate, below of=u_end, node distance=(\hLTI+\dyLTI)]{};
	\node (front_LTI1_c)	[coordinate, right of=u1_start, node distance=(\wShortNamedSignal+\wLTI/2)]{};
	\node (front_LTI2_c)	[coordinate, right of=u_end, node distance=(\wShortNamedSignal+\wLTI/2)]{};
	\node (front_LTI3_c)	[coordinate, right of=u3_start, node distance=(\wShortNamedSignal+\wLTI/2)]{};
	\node (SNL1_c)	[coordinate, right of=front_LTI1_c, node distance=(\wLTI/2+\wShortNamedSignal+\wSNL/2)]{};
	\node (SNL2_c)	[coordinate, right of=front_LTI2_c, node distance=(\wLTI/2+\wShortNamedSignal+\wSNL/2)]{};
	\node (SNL3_c)	[coordinate, right of=front_LTI3_c, node distance=(\wLTI/2+\wShortNamedSignal+\wSNL/2)]{};
	\node (back_LTI1_c)	[coordinate, right of=SNL1_c, node distance=(\wSNL/2+\wShortNamedSignal+\wLTI/2)]{};
	\node (back_LTI2_c)	[coordinate, right of=SNL2_c, node distance=(\wSNL/2+\wShortNamedSignal+\wLTI/2)]{};
	\node (back_LTI3_c)	[coordinate, right of=SNL3_c, node distance=(\wSNL/2+\wShortNamedSignal+\wLTI/2)]{};
	\node (y1_end)	[coordinate, right of=back_LTI1_c, node distance=(\wLTI/2+\wShortNamedSignal)]{};
	\node (y2_end)	[coordinate, right of=back_LTI2_c, node distance=(\wLTI/2+\wShortNamedSignal)]{};
	\node (y3_end)	[coordinate, right of=back_LTI3_c, node distance=(\wLTI/2+\wShortNamedSignal)]{};
	\node (sum_c)		[coordinate, right of=y2_end, node distance=(\wUnnamedSignal+\sSum/2)]{};
	\node (y_end)		[coordinate, right of=sum_c, node distance=(\sSum/2+\wShortNamedSignal)]{};

	% specific nodes
        \node (front_LTI1)	at (front_LTI1_c) 	[SISO_LTI]{$P_1(z)$};
	\node (front_LTI2)	at (front_LTI2_c) 	[SISO_LTI, draw=none, black]{};
        \node (front_LTI3)	at (front_LTI3_c) 	[SISO_LTI]{$P_2(z)$};
        \node (SNL1)		at (SNL1_c)		[SISO_BigSNL]{$f(x)=x^3$};
	\node (SNL2)		at (SNL2_c)		[SISO_BigSNL, draw=none, black]{};
        \node (SNL3)		at (SNL3_c)		[SISO_BigSNL]{$f(x)=x^3$};
        \node (back_LTI1)	at (back_LTI1_c) 	[SISO_LTI]{$Q_1(z)$};
	\node (back_LTI2)	at (back_LTI2_c) 	[SISO_LTI, draw=none, black]{};
        \node (back_LTI3)	at (back_LTI3_c) 	[SISO_LTI]{$Q_2(z)$};
	\node (sum)		at (sum_c)		[sum]{};
	\node 			at (sum)		{$+$};

	% signals
	\path[-]	(u_start)	-- node {$u(t)$}	(u_end);
	\path[-]	(sum)		-- node {$y(t)$}	(y_end);

	% arrows
	\draw[->] 	(u_start)	-- 	(u_end);
	\draw[<->] 	(front_LTI1)	-- 	(u1_start)	-- 	(u3_start)	-- 	(front_LTI3);
	\draw[->] 	(front_LTI1)	-- 	(SNL1);
	\draw[->] 	(front_LTI3)	-- 	(SNL3);
	\draw[->] 	(SNL1)		-- 	(back_LTI1);
	\draw[->] 	(SNL3)		-- 	(back_LTI3);
	\draw[->] 	(back_LTI1)	-- 	(y1_end)	-- 	(sum);
	\draw[->] 	(back_LTI3)	-- 	(y3_end)	-- 	(sum);
	\draw[->] 	(sum)		-- 	(y_end);
\end{tikzpicture}

%% file: fig-pWHnoise.tex
% parallel wiener-hammerstein figure 
\begin{tikzpicture}[font=\Large, >=stealth', thick, auto]
	% general coordinates
	\node (u_start)	[coordinate]{};
	\node (u_end)		[coordinate, right of=u_start, node distance=\wShortNamedSignal]{};
	\node (u1_start)	[coordinate, above of=u_end, node distance=(\hLTI+\dyLTI)]{};
	\node (u3_start)	[coordinate, below of=u_end, node distance=(\hLTI+\dyLTI)]{};
	\node (front_LTI1_c)	[coordinate, right of=u1_start, node distance=(\wShortNamedSignal+\wLTI/2)]{};
	\node (front_LTI2_c)	[coordinate, right of=u_end, node distance=(\wShortNamedSignal+\wLTI/2)]{};
	\node (front_LTI3_c)	[coordinate, right of=u3_start, node distance=(\wShortNamedSignal+\wLTI/2)]{};
	\node (SNL1_c)	[coordinate, right of=front_LTI1_c, node distance=(\wLTI/2+\wShortNamedSignal+\wSNL/2)]{};
	\node (SNL2_c)	[coordinate, right of=front_LTI2_c, node distance=(\wLTI/2+\wShortNamedSignal+\wSNL/2)]{};
	\node (SNL3_c)	[coordinate, right of=front_LTI3_c, node distance=(\wLTI/2+\wShortNamedSignal+\wSNL/2)]{};
	\node (back_LTI1_c)	[coordinate, right of=SNL1_c, node distance=(\wSNL/2+\wShortNamedSignal+\wLTI/2)]{};
	\node (back_LTI2_c)	[coordinate, right of=SNL2_c, node distance=(\wSNL/2+\wShortNamedSignal+\wLTI/2)]{};
	\node (back_LTI3_c)	[coordinate, right of=SNL3_c, node distance=(\wSNL/2+\wShortNamedSignal+\wLTI/2)]{};
	\node (y1_end)	[coordinate, right of=back_LTI1_c, node distance=(\wLTI/2+\wShortNamedSignal)]{};
	\node (y2_end)	[coordinate, right of=back_LTI2_c, node distance=(\wLTI/2+\wShortNamedSignal)]{};
	\node (y3_end)	[coordinate, right of=back_LTI3_c, node distance=(\wLTI/2+\wShortNamedSignal)]{};
	\node (sum_c)		[coordinate, right of=y2_end, node distance=(\wUnnamedSignal+\sSum/2)]{};
	\node (y_end)		[coordinate, right of=sum_c, node distance=(\sSum/2+\wShortNamedSignal)]{};

        \node (sum_end)         [coordinate, right of=sum_c, node distance=(\wShortNamedSignal+\sSum)]{};
	\node (yn_end)		[coordinate, right of=sum_end, node distance=(\sSum/2+\wShortNamedSignal)]{};
        \node (ny_start)        [coordinate, above of=sum_end, node distance=(1.2*\wShortNamedSignal)]{};

	% specific nodes
        \node (front_LTI1)	at (front_LTI1_c) 	[SISO_LTI]{$P_1(z)$};
	\node (front_LTI2)	at (front_LTI2_c) 	[SISO_LTI, draw=none, black]{};
        \node (front_LTI3)	at (front_LTI3_c) 	[SISO_LTI]{$P_2(z)$};
        \node (SNL1)		at (SNL1_c)		[SISO_SNL]{$f_1(\cdot)$};
	\node (SNL2)		at (SNL2_c)		[SISO_SNL, draw=none, black]{};
        \node (SNL3)		at (SNL3_c)		[SISO_SNL]{$f_2(\cdot)$};
        \node (back_LTI1)	at (back_LTI1_c) 	[SISO_LTI]{$Q_1(z)$};
	\node (back_LTI2)	at (back_LTI2_c) 	[SISO_LTI, draw=none, black]{};
        \node (back_LTI3)	at (back_LTI3_c) 	[SISO_LTI]{$Q_2(z)$};
	\node (sum)		at (sum_c)		[sum]{};
	\node 			at (sum)		{$+$};
        \node (sum2)            at (sum_end)            [sum]{};
        \node                   at (sum2)               {$+$};
        \node[anchor=north west] at (ny_start)          {$n_y(k)$}; 

	% signals
	\path[-]	(u_start)	-- node {$u(t)$}	(u_end);
	\path[-]	(sum)		-- node {$y_0(t)$}	(y_end);
        \path[-]        (sum2)          -- node {$y(k)$}        (yn_end);

	% arrows
	\draw[->] 	(u_start)	-- 	(u_end);
	\draw[<->] 	(front_LTI1)	-- 	(u1_start)	-- 	(u3_start)	-- 	(front_LTI3);
	\draw[->] 	(front_LTI1)	-- 	(SNL1);
	\draw[->] 	(front_LTI3)	-- 	(SNL3);
	\draw[->] 	(SNL1)		-- 	(back_LTI1);
	\draw[->] 	(SNL3)		-- 	(back_LTI3);
	\draw[->] 	(back_LTI1)	-- 	(y1_end)	-- 	(sum);
	\draw[->] 	(back_LTI3)	-- 	(y3_end)	-- 	(sum);
	\draw[->] 	(sum)		-- 	(sum2);
        
        \draw[->]       (sum2)          --      (yn_end);
        \draw[->]       (ny_start)      --      (sum2);
\end{tikzpicture}